\documentclass[12pt]{article}

\usepackage{graphicx,amssymb}

\newtheorem{theorem}{Theorem}[section]
\newtheorem{lemma}[theorem]{Lemma}

\newcommand{\bitem}{\begin{itemize}}
\newcommand{\eitem}{\end{itemize}}
\newcommand{\bcenter}{\begin{center}}
\newcommand{\ecenter}{\end{center}}
\newcommand{\benum}{\begin{enumerate}}
\newcommand{\eenum}{\end{enumerate}}
\newcommand{\beqar}{\begin{eqnarray*}}
\newcommand{\eeqar}{\end{eqnarray*}}
\newcommand{\beq}{\begin{equation}}
\newcommand{\eeq}{\end{equation}}
\newcommand{\bdesc}{\begin{description}}
\newcommand{\edesc}{\end{description}}

\newcommand{\Bin}{{\rm Bin}}

\newcommand{\eps}{\varepsilon}
\newcommand{\graph}{{\rm graph}}
\newcommand{\goto}{\rightarrow}

\newcommand{\Hold}{{\rm H}}

\newcommand{\cF}{{\cal F}}

\newcommand{\cH}{{\cal H}}

\newcommand{\cM}{{\cal M}}

\newcommand{\bbS}{{\mathbb{S}}}

\newcommand{\bS}{{\bf S}}
\newcommand{\bs}{{\bf s}}
\newcommand{\bt}{{\bf t}}
\newcommand{\bn}{{\bf n}}
\newcommand{\bm}{{\bf m}}
\newcommand{\bzero}{{\bf 0}}

\newcommand{\bbO}{\mathbb{O}}
\newcommand{\bbG}{\mathbb{G}}
\newcommand{\bbN}{\mathbb{N}}

\newcommand{\bbR}{\mathbb{R}}
\newcommand{\bbB}{\mathbb{B}}

\newcommand{\inner}[2]{\langle #1,#2 \rangle}
\newcommand{\ang}[2]{\angle \left(#1,#2\right)}
\newcommand{\Span}[1]{{\rm span}\left\{#1\right\}}

\newcommand{\pr}[1]{{\bf P}\left\{#1\right\}}

\newcommand{\qed}{{\unskip\nobreak\hfil\penalty50\hskip2em\vadjust{}
\nobreak\hfil$\Box$\parfillskip=0pt\finalhyphendemerits=0\par}\vspace{.1cm}}

\newcommand{\expect}[1]{{\bf E}\left\{#1\right\}}
\newcommand{\var}[1]{{\bf var}\left\{#1\right\}}

\newcommand{\acos}{{\rm acos}}

\newcommand{\bino}[2]{\left(\begin{array}{c} #1 \\ #2
    \end{array}\right)}

\begin{document}




\title{Interpolation of Random Hyperplanes}
\author{Ery Arias-Castro
\\ 
University of California, San Diego}
\date{}
\maketitle

\begin{abstract}
Let $\{(Z_i,W_i):i=1,\dots,n\}$ be uniformly distributed in $[0,1]^d
\times \bbG(k,d)$, 
where $\bbG(k,d)$ denotes the space of
$k$-dimensional linear subspaces of $\bbR^d$.
For a differentiable function $f: [0,1]^k \goto [0,1]^d$, we say that $f$ 
interpolates $(z,w) \in [0,1]^d \times \bbG(k,d)$ if there exists $x
\in [0,1]^k$ such that $f(x) = z$ and $\vec{f}(x) = w$, where
$\vec{f}(x)$ denotes the tangent space at $x$ defined by $f$.
For a smoothness class $\cF$ of H\"older type, 
we obtain probability bounds on the maximum number of points a
function $f \in \cF$ interpolates.
\end{abstract}

\footnotetext[1]{This work was partially supported by NSF grant DMS-0603890.  The author was at the Mathematical Sciences Research Institute while preparing the first draft.}
\footnotetext[2]{
The author would like to thank Emmanuel Cand\`es, David Donoho, Bruce Driver, Bo'az Klartag and Allen Knutson for helpful discussions and references.}




\section{Introduction}

\begin{figure}[htbp]
\centering
\begin{tabular}{cc}
(a) Under $H_0$ & (b) Under $H_1$\\
\includegraphics[height=2in]{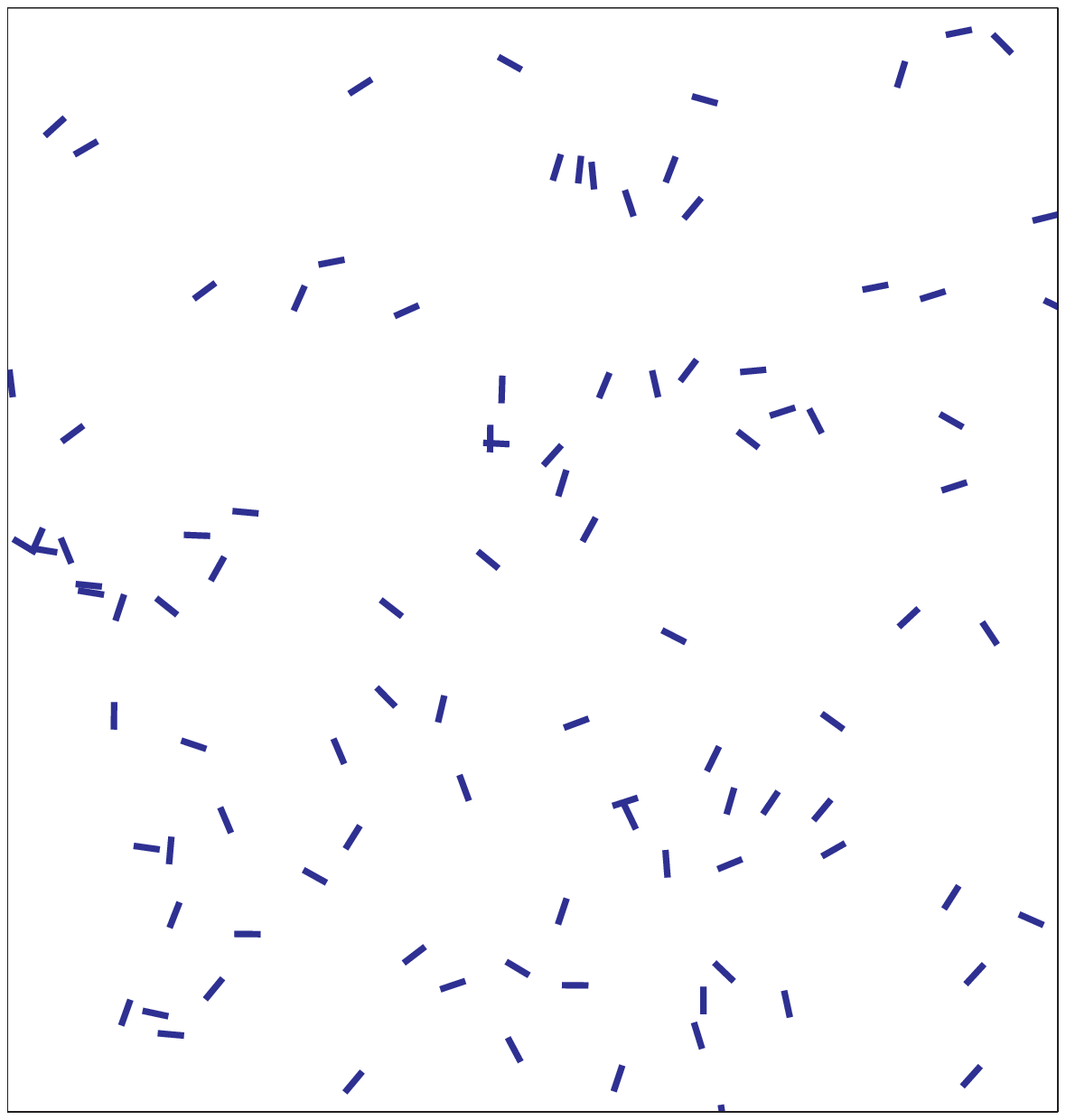} &
\includegraphics[height=2in]{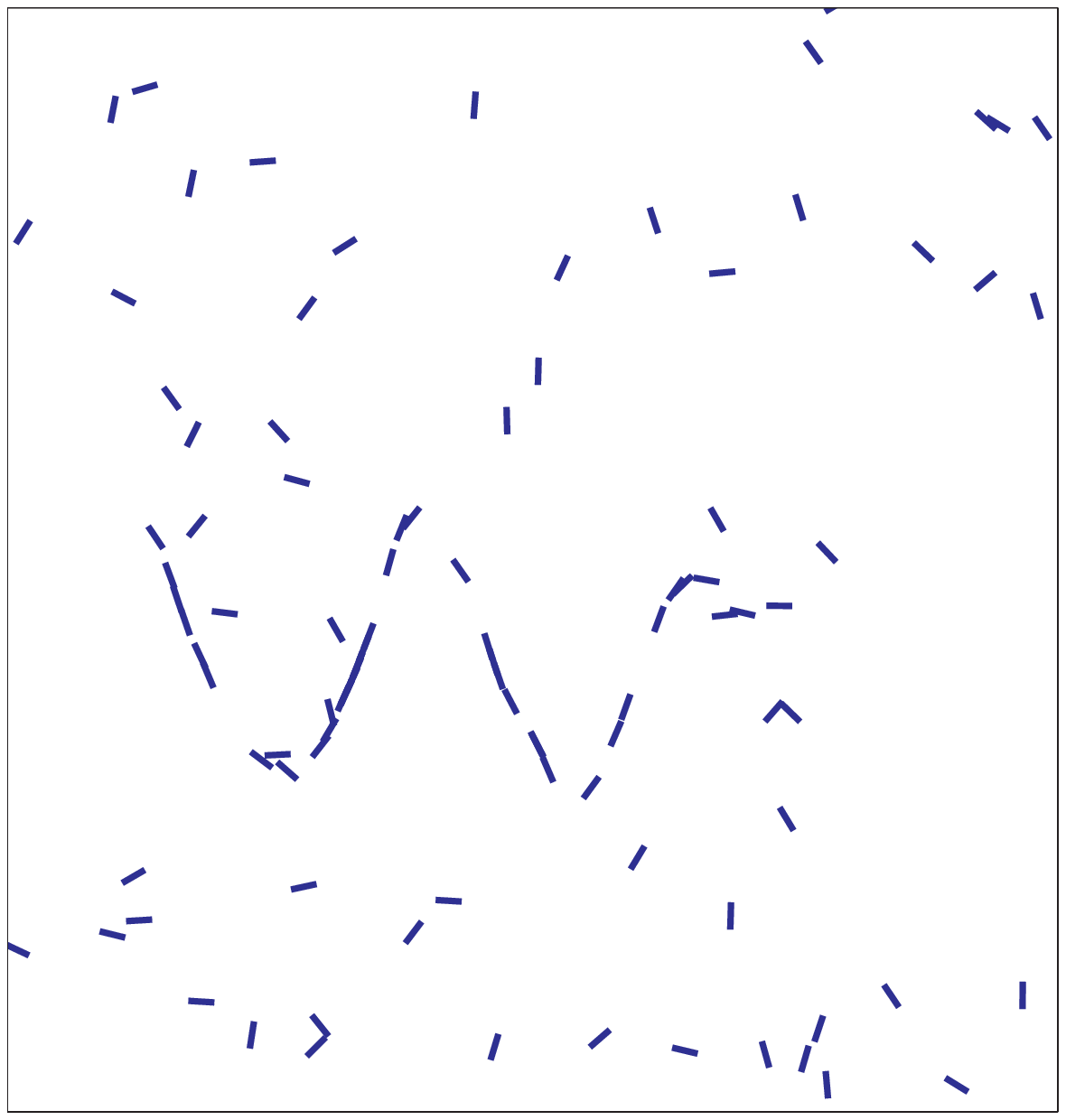}
\end{tabular}
\caption{In Panel (a) we observe a realization under the null
  hypothesis ($n = 100$). In Panel (b) we observe a realization under
  the alternative hypothesis ($n = 100,n_1 = 40$).} 
\label{fig:direction}
\end{figure}

This paper is motivated by experiments in the field of Psychophysics
\cite{FieHayHes} that study the ability of the Human Visual System at
detecting curvilinear features in background clutter. 
In these experiments, human subjects are shown an image consisting of
oriented small segments of same length dispersed in a square, such as
in Figure \ref{fig:direction}.  

The locations and orientations of these segments are either purely random
(panel (a)) or a curve is actually ``hidden'' among purely random clutter,
which here means that a curve was used to simulate a fraction of the
segments by randomly choosing segments that are tangent to the curve
at their midpoint (panel (b)).

From a Statistics viewpoint, 
this detection task, that human subjects are asked to perform, can be
formalized into a hypothesis testing problem.

We say that a curve $\gamma \subset [0,1]^2$, parametrized by arclength, interpolates $(z,w) \in
[0,1]^2 \times \bbS^1$ if there is $x$ such that $\gamma(x) = z$ and
$\dot{\gamma}(x) = w$, where $\bbS^1$ denotes the unit circle and $\dot{\gamma}(x)$ the derivative of $\gamma$ at $x$.\\

We observe $n$ segments of fixed length dispersed in
the unit square. 
\bitem
\item
Under the {\it null hypothesis}, the segments have locations and
orientations sampled uniformly at random in $[0,1]^2 \times \bbS^1$.  
\item
Under the {\it (composite) alternative hypothesis}, the
segments are as above except for $n_1$ of them that are chosen among 
those that a fixed curve $\gamma$
interpolates.   
The curve $\gamma$ is unknown but restricted to
belong to some known class $\Gamma$.
\eitem  
Note that we do not specify the distribution of the segments tangent
to the curve.\\
  
For $\gamma \in \Gamma$, define
$$N_n^\goto(\gamma) = \#\{i=1,\dots,n: \gamma\ {\rm interpolates}\
(Z_i,W_i)\},$$
and, with some abuse of notation,
$$N_n^\goto(\Gamma) = \max_{\gamma \in \Gamma} N_n^\goto(\gamma).$$

In \cite{AriDonHuoTov}, the test that rejects for large
$N_n^\goto(\Gamma)$ was analyzed for $\Gamma$ the class of curves in
the unit square with length and curvature bounded by some constant
$c > 0$. 
In particular, it was shown that, under the null hypothesis, 
for some constants $0 < A < B <
\infty$, 
$$\pr{A\ n^{1/4} \leq N_n^\goto(\Gamma) \leq B\ n^{1/4}} \goto 1,
\quad n \goto \infty.$$
Note that the upper bound implies that this test is powerful when $n_1 \geq B n^{1/4}$.\\ 



In this paper, we generalize this setting to higher dimensions.
Let $\bbG(k,d)$ be the set of $k$-dimensional linear subspaces in $\bbR^d$. 
To $\bbG(k,d)$ we associate its uniform measure $\lambda$, which is the only invariant probability measure on $\bbG(k,d)$ that is invariant under the action of the orthogonal group $\bbO(d)$ -- see \cite{MilSch}, Section 1.  

For a function $f:[0,1]^k \goto [0,1]^d$ differentiable at $x$, let
$$\vec{f}(x) = \Span{\partial_s f(x): s=1,\dots,k}.$$

A function $f:[0,1]^k \goto [0,1]^d$ is said to interpolate $(z,w) \in
[0,1]^d \times \bbG(k,d)$ if there exists $x \in [0,1]^k$ such that
$f(x) = z$ and $\vec{f}(x) = w$.\\

We consider the following hypothesis testing problem. 
We observe 
$$\{(Z_i,W_i): i=1,\dots,n\} \subset [0,1]^d \times
\bbG(k,d).$$ 

\bitem
\item
Under the {\it null hypothesis}, $\{(Z_i,W_i): i=1,\dots,n\}$ are 
independent and identically uniformly distributed in $[0,1]^d \times
\bbG(k,d)$.
\item
Under the {\it (composite) alternative hypothesis}, $\{(Z_i,W_i):
i=1,\dots,n\}$ are as above except for $n_1$ of them that are chosen 
among those that a fixed function $f$
interpolates.
The function $f$ is unknown but restricted to
belong to some known class $\cF$. 
\eitem

Before specifying $\cF$, we introduce some notation.
For a vector $x = (x_1,\dots,x_d) \in \bbR^d$, the supnorm is defined as $\|x\|_\infty = \max\{|x_i|:i=1,\dots,d\}$.
For a function $f:\Omega \subset \bbR^k \goto \bbR^d$, $\|f\|_\infty = \sup_{x \in \Omega} \|f(x)\|_\infty$.
The Euclidean inner product and the corresponding norm are denoted by $\inner{\cdot}{\cdot}$ and $\|\cdot\|$ respectively.
The angle $\ang{H}{K} \in [0,\pi]$ between two linear subspaces $H,K \subset \bbR^d$, with $1 \leq {\rm dim} H \leq {\rm dim} K$, is defined by
$$\ang{H}{K} = \max_{u \in H} \min_{v \in K} \acos\left(\frac{\inner{u}{v}}{\|u\| \|v\|}\right).$$
This corresponds to the largest canonical angle as defined in \cite{GolVan} and constitutes a metric on $\bbG(k,d)$ -- see also \cite{AbsEdeKoe} for a related study of the largest canonical angle between two subspaces uniformly distributed in $\bbG(k,d)$.\\

The class $\cF$, parametrized by $\beta \geq 1$, is defined as the set of twice differentiable,
one-to-one functions $f:[0,1]^k \goto [0,1]^d$ with the following
additional properties:
\bitem
\item
For all $s = 1,\dots,k$, $1/\beta \leq \|\partial_s f (x)\|_\infty \leq \beta$ for all $x \in [0,1]^k$;    
\item 
For all $s,t = 1,\dots,k$, $\|\partial_{st} f\|_\infty \leq \beta$;
\item
For all $s = 1,\dots,k$ and $x \in [0,1]^k$,
$$\ang{\partial_s f(x)}{\Span{\partial_t f(x): t \neq s}}
\geq \frac{1}{2 \beta (d-k)},$$ 
which is void if $k = 1$.
(In this paper, we identify a non-zero vector with the one dimensional linear subspace it generates.)
\eitem
The last condition and the constraint $\beta \geq 1$ ensure that $\cF$ contains graphs of the form $x \goto (x,g(x))$, where $g:[0,1]^k \goto [0,1]^{d-k}$ satisfies the first two conditions -- see Lemma \ref{lem:HinF}.

Define
$$N_n^{\goto}(f) = \#\ \{i=1,\dots,n: f\ {\rm interpolates}\
(Z_i,W_i)\},$$
and, with some abuse of notation,
$$N_n^{\goto}(\cF) = \max_{f \in \cF}\ N_n^{\goto}(f).$$
  
Let 
$$\vec{\rho} = \frac{k}{k+(d-k)(k+2)}.$$

\begin{theorem}
\label{th:direction_glrt_ub}
There is a constant $B = B(k,d,\beta) < \infty$ such that, under the null
hypothesis, 
$$\pr{N_n^{\goto}(\cF) > B\ n^{\vec{\rho}}} \goto 0,
\quad n \goto \infty.$$ 
\end{theorem}
As before, this implies that the test that rejects for large
values of $N_n^{\goto}(\cF)$ is powerful when $n_1 > B n^{\vec{\rho}}$.

\begin{theorem}
\label{th:direction_glrt_lb}
There is a constant $A  = A(k,d,\beta) > 0$ such that, under the null hypothesis,  
$$\pr{N_n^{\goto}(\cF) < A\ n^{\vec{\rho}}} \goto 0, \quad n \goto
\infty.$$   
\end{theorem}

The remaining of the paper is organized as follows.  
In Section 2 we introduce a related, yet different hypothesis testing
problem.   
In Section 3 and Section 4, we prove results announced in Section 2.  
In Section 5 and Section 6, we follow the arguments in Section 3 and Section 4 to prove Theorem
\ref{th:direction_glrt_ub} and Theorem 
\ref{th:direction_glrt_lb}.  
Some intermediary lemmas are proved in the Appendix.

\section{Another Hypotheses Testing Problem}

We introduce another hypothesis testing problem as a stepping stone
towards proving Theorem \ref{th:direction_glrt_ub} and Theorem
\ref{th:direction_glrt_lb}, and also for its own sake.\\

Let $\alpha > 1$, $\beta>0$ and define $r = \lfloor \alpha
\rfloor = \max\{m \in \bbN: m < \alpha\}$.  (In this paper, we include $0$ in $\bbN$.)
Define the H\"older class $\Hold^{k,d}(\alpha,\beta)$ to be the set
of functions $f:[0,1]^k 
\goto [0,1]^d$, with $f = (f_1,\dots,f_d)$ such that,
for all $\bs = (s_1,\dots,s_k) \in \bbN^k$ with $|\bs| = s_1 +
\cdots + s_k \leq r$,
$$\|f^{(\bs)}\|_\infty \leq \beta;$$ 
and, for all $\bs \in \bbN^k$ with $|\bs| = r$, 
$$\|f^{(\bs)}(x) - f^{(\bs)}(y)\|_\infty \leq \beta \|x -
y\|_\infty^{\alpha-r},$$
where $f^{(\bs)} = (\partial_{s_1 \cdots s_k} f_1, \dots,
\partial_{s_1 \cdots s_k} f_d)$.
When there is no possible confusion, we use the notation $\cH =
\Hold^{k,d-k}(\alpha,\beta)$.\\ 

Fix $r_0$ an integer such that
$1 \leq r_0 \leq r$.  
Let 
$$\bS = \{\bs \in \bbN^k: |\bs| \leq r_0\},$$
with cardinality
$$|\bS| = \sum_{s=0}^{r_0} \bino{s+k-1}{k-1}.$$

We denote by $y^{\bS}$ a vector in $\bbR^{(d-k)|\bS|}$ and by
$f^{(\bS)}(x)$ the vector $(f^{(\bs)}(x): \bs \in \bS)$.
A function $f \in \cH$ is said to interpolate
$(x,y^{\bS}) \in [0,1]^k \times \bbR^{(d-k)|\bS|}$ if $f^{(\bS)}(x) =
y^{\bS}$.\\

Consider the following hypothesis testing problem.
We observe 
$$\{(X_i,Y_i^{\bS})): i=1,\dots,n\} \subset
[0,1]^k \times \bbR^{(d-k)|\bS|}.$$

\begin{figure}[htbp]
\centering
\begin{tabular}{cc}
(a) Under $H_0$ & (b) Under $H_1$\\
\includegraphics[height=2in]{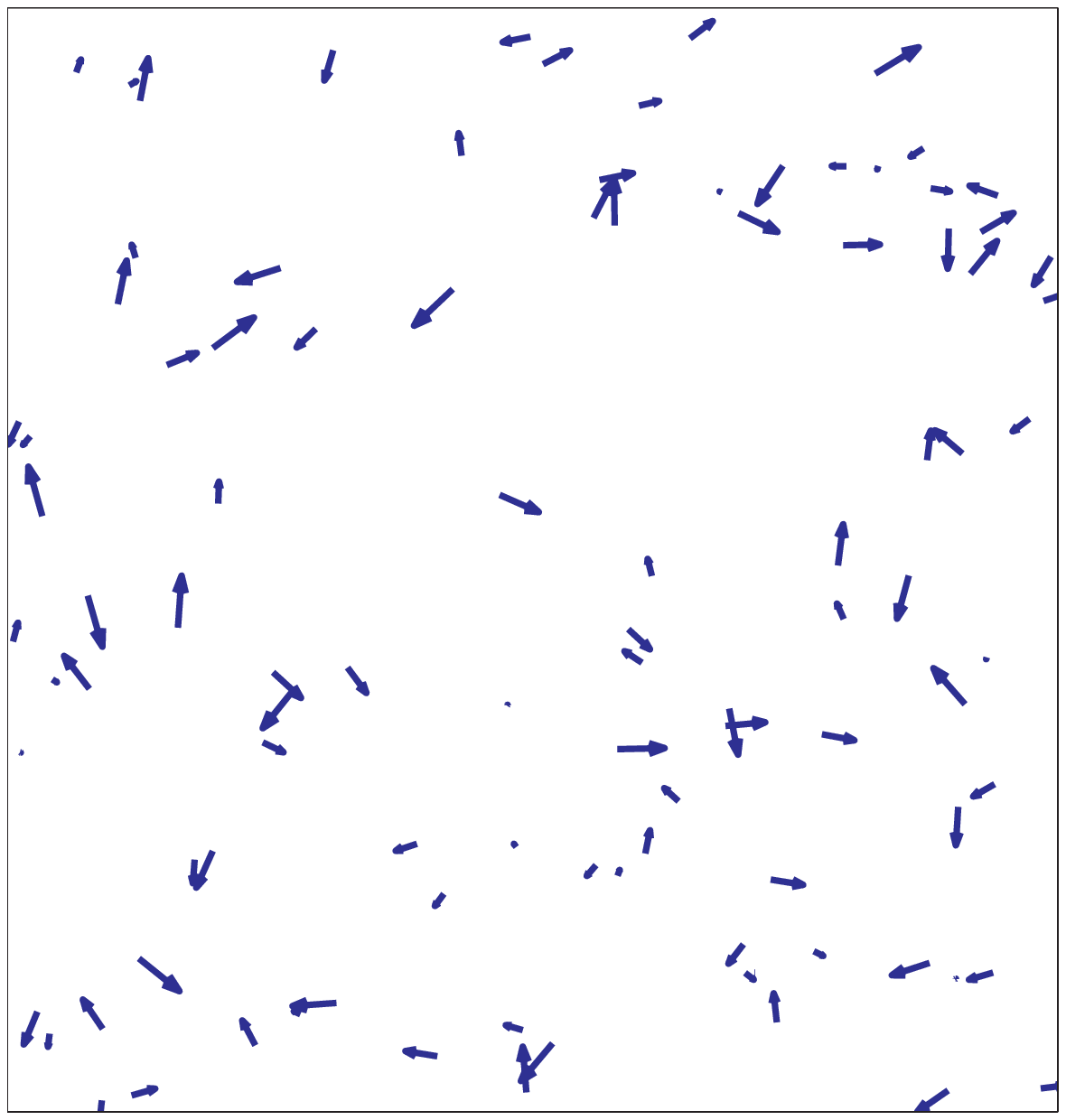} &
\includegraphics[height=2in]{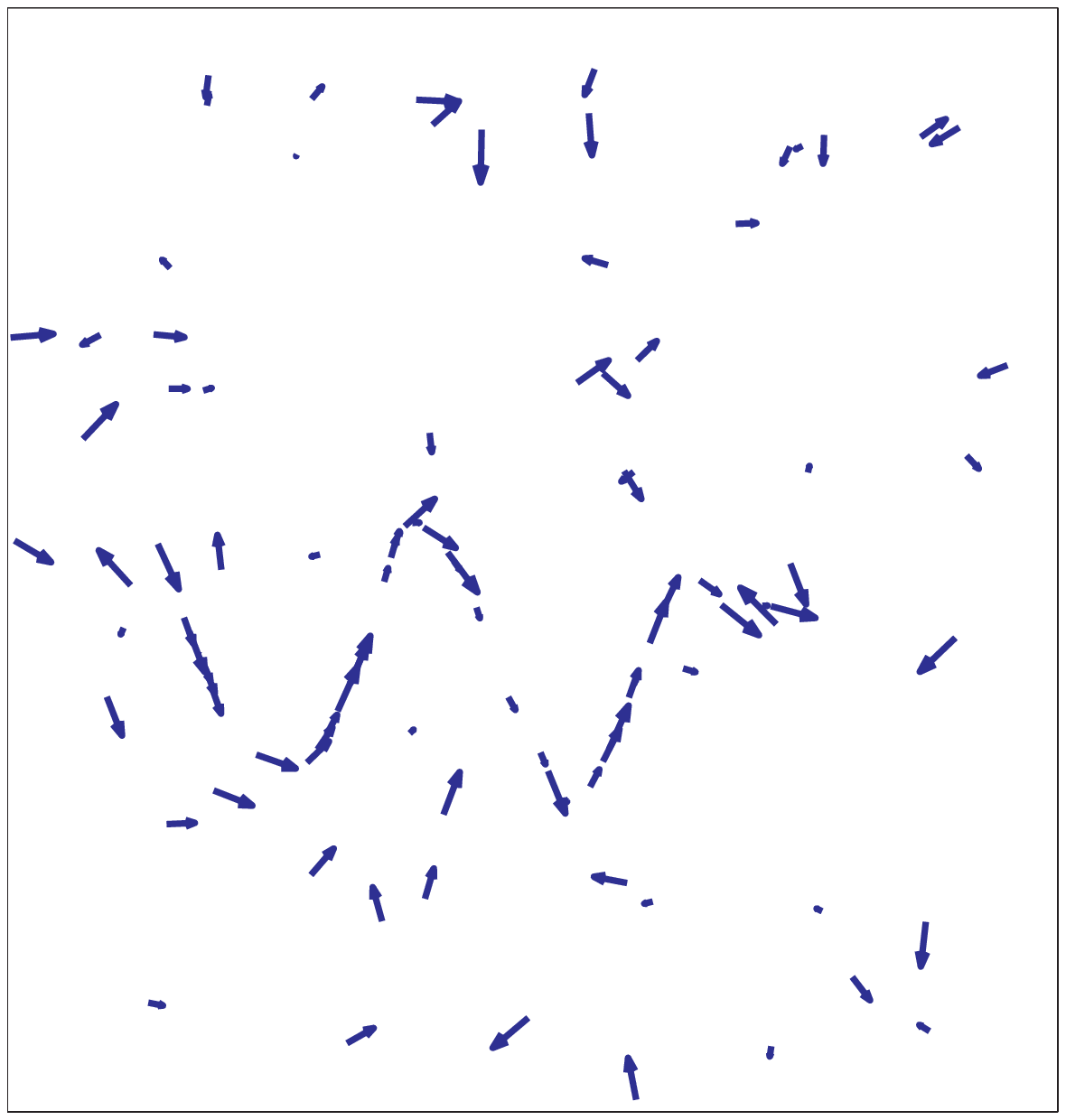}
\end{tabular}
\caption{In Panel (a) we observe a realization under the null
  hypothesis ($n = 100$). In Panel (b) we observe a realization under
  the alternative hypothesis ($n = 100,n_1 = 40$).} 
\label{fig:velocity}
\end{figure}

\bitem
\item Under the {\it null hypothesis}, $\{(X_i,Y_i^{\bS}):
  i=1,\dots,n\}$ are independent and identically uniformly distributed
  in $[0,1]^k \times [0,1]^{d-k} \times
  [-\beta,\beta]^{(d-k)(|\bS|-1)}$. 
\item Under the {\it (composite) alternative hypothesis},
  $\{(X_i,Y_i^{\bS}): i=1,\dots,n\}$ are as above
  except for $n_1$ of them that are chosen among those that a fixed
  function $f$ interpolates.
  The function $f$ is unknown but restricted to
  belong to $\cH$.
\eitem
Figure \ref{fig:velocity} shows an example, with $d = 2$ and $r_0 = 1$.

Define
$$N_n^{(r_0)}(f) = \# \{i=1,\dots,n: f\ {\rm interpolates}\
(X_i,Y_i^{\bS})\},$$
and, with some abuse of notation,
$$N_n^{(r_0)}(\cH) = \max_{f \in \cH} N^{(r_0)}(f).$$ 

Let 
$$\rho(r_0) = \frac{k}{k + \alpha (d - k) w} =
\frac{1}{1 + \alpha (d/k - 1) w}\ ,$$
where
$$w = \sum_{s=0}^{r_0} (1-s/\alpha) \bino{s+k-1}{k-1}.$$ 

\begin{theorem}
\label{th:velocity_glrt_ub}
There is a constant $B = B(k,d,\alpha,\beta,r_0) > 0$ such that, under the
null hypothesis,
$$\pr{N_n^{(r_0)}(\cH) > B\ n^{\rho(r_0)}} \goto 0,
\quad n \goto \infty.$$ 
\end{theorem}
As before, this implies that the test that rejects for large
values of $N_n^{(r_0)}(\cH)$ is powerful when $n_1 > B n^{\rho(r_0)}$.

\begin{theorem}
\label{th:velocity_glrt_lb}
There is a constant $A  = A(k,d,\alpha,\beta,r_0) > 0$ such that, under the null hypothesis,
$$\pr{N_n^{(r_0)}(\cH) < A\ n^{\rho(r_0)}} \goto 0, \quad n \goto
\infty.$$   
\end{theorem}

\noindent {\bf Remark.} 
For $\alpha = 2$ and $r_0=1$, $\rho(r_0) = \vec{\rho}$, meaning that $N_n^{(r_0)}(\cH)$ and $N_n^\goto(\cF)$ are, in that case, of same order of magnitude with high probability. 
This will be used explicitly in Section 6 when proving Theorem \ref{th:direction_glrt_lb}.

\section{Proof of Theorem \ref{th:velocity_glrt_ub}}
\label{proof:velocity_glrt_ub}

For $y_1^\bS,y_2^\bS \in \bbR^{(d-k)|\bS|}$, define the
discrepancy 
$$\Phi(y_1^\bS,y_2^\bS) = \max_{\bs \in \bS} \|y_1^\bs -
y_2^\bs\|_\infty^{\alpha/(\alpha-|\bs|)}.$$ 

The discrepancy $\Phi$ induces a discrepancy on functions, namely
$$\Phi(f,g) = \max_{\bs \in \bS}
\|f^{(\bs)}-g^{(\bs)}\|_\infty^{\alpha/(\alpha-|\bs|)}.$$ 

The argument for proving Theorem \ref{th:velocity_glrt_ub} is based on
coverings of $\cH$ with respect to $\Phi$.
Let $L_\eps$ be the $\eps$-covering number of $\cH$ with respect to $\Phi$.

\begin{lemma}
\label{lem:entropy}
There is a constant $c_1 = c_1(k,d,\alpha,\beta,r_0) > 0$ such that, for all $\eps > 0$,
$$\log L_\eps \leq c_1 \eps^{-k/\alpha}.$$
\end{lemma}
Lemma \ref{lem:entropy} follows immediately from the proof of Theorem XIII in
\cite{Kol}, Chapter ``$\eps$-entropy and $\eps$-capacity of sets in
functional spaces''.\\ 

For a set $K \subset \bbR^k \times \bbR^{(d-k)|\bS|}$ and $\eps > 0$,
we denote by $K^{\eps,\Phi}$ the set of points $(x,y^\bS)$ such that
there is $(x_1,y_1^\bS) \in K$ with $x_1 = x$ and $\Phi(y_1^\bS,y^\bS)
\leq \eps$. 

For each $\eps$, we select an $\eps$-net
$\{f_j:j=1,\dots,L_\eps\}$ of $\cH$ for $\Phi$.
For $j=1,\dots,L_\eps$, we define
\beqar
K_j 
& = & \{(x,y^{\bS}) \in \bbR^k \times \bbR^{(d-k)|\bS|}: \Phi(y^{\bS},f_j^{(\bS)}(x)) \leq \eps\}\\
& = & \graph^{\bS}(f_j)^{\eps,\Phi},
\eeqar
where, for $f \in \cH$,
$$\graph^{\bS}(f) = \{(x,f^{(\bS)}(x)): x \in [0,1]^k\} \subset [0,1]^k
\times \bbR^{(d-k)|\bS|}.$$

We extend $N^{(r_0)}(\cdot)$ to subsets $K \subset \bbR^k \times \bbR^{(d-k)|\bS|}$, by defining $N^{(r_0)}(K)$ to be the number of points $(X_i,Y_i^{\bS})$ that belong to $K$.

Let
$$M_n(\eps) = \max_{j=1,\dots,L_\eps}\
N_n^{(r_0)}(K_j).$$  
By definition, it is straightforward to see that 
$$N_n^{(r_0)}(\cH) \leq M_n(\eps),$$ 
for all $\eps > 0$.
We therefore focus on bounding $M_n(\eps)$.\\   

By Boole's inequality, we have
$$\pr{M_n(\eps) > b} \leq L_\eps \cdot
\max_{j=1,\dots,L_\eps}\ \pr{N_n^{(r_0)}(K_j) > b}.$$ 
Moreover, we know that, for any set $K \subset \bbR^k \times
\bbR^{(d-k)|\bS|}$,  
$$N_n^{(r_0)}(K) \backsim \Bin(n,\mu(K)),$$ 
where $\mu$ is the uniform
measure on $[0,1]^k \times [0,1]^{d-k} \times
[-\beta,\beta]^{(d-k)(|\bS|-1)}$. 

Hence, 
\beqar
\pr{M_n(\eps) > b} 
& \leq & L_\eps \cdot \max_{j=1,\dots,L_\eps}\
\pr{\Bin(n,\mu(K_j)) > b}\\
& = & L_\eps \cdot \pr{\Bin(n,\max_j \mu(K_j)) > b}.  
\eeqar 

\begin{lemma}
\label{lem:graphS_measure}
There is $c_2 > 0$ such that, for all $f \in \cH$ and all $\eps > 0$, 
$$\mu\left(\graph^{\bS}(f)^{\eps,\Phi}\right) \leq c_2
\eps^{(d-k)w}.$$  
\end{lemma}
{\it Proof of Lemma \ref{lem:graphS_measure}.}
Start with 
$$\graph^{\bS}(f)^{\eps,\Phi} \subset \cup_{x \in [0,1]^k}\ \{x\}
\bigotimes_{\bs \in \bS} \bbB(f^{(\bs)}(x),\eps^{1 - |\bs|/\alpha}),$$
where $\bbB(y,\eta)$ denotes the supnorm ball in $\bbR^{d-k}$ centered
at $y$ of radius $\eta$.

Hence, integrating over $x \in [0,1]^k$ last, we have
$$\mu(\graph^{\bS}(f)^{\eps,\Phi}) \leq c_2 \prod_{\bs \in \bS}\ (\eps^{1 -
  |\bs|/\alpha})^{d-k} = c_2\ \eps^{(d-k)w}.$$
\qed   

Using Lemma \ref{lem:graphS_measure}, we arrive at
$$\pr{M_n(\eps) > b} \leq L_\eps \cdot \pr{\Bin(n,c_2
  \eps^{(d-k)w}) > b}.$$   

\begin{lemma}
\label{lem:bin_tail}
There is a constant $c > 0$ such that, for any $n$ positive integer
and $p \in (0,1/2)$, and for all $b > 2 np$, 
$$\pr{\Bin(n,p) > b} \leq \exp(-c\cdot b).$$  
\end{lemma}
Lemma \ref{lem:bin_tail} follows directly from standard large
deviations bounds for binomial random variables -- see \cite{ShoWel},
p. 440, Inequality 1, (ii). 

We use Lemma \ref{lem:bin_tail} to obtain
$$\pr{\Bin(n,c_2 \eps^{(d-k)w}) > b} \leq \exp(- c\cdot b),$$
for all $b > 2 n c_2 \eps^{(d-k)w}$. 

Collecting terms, we arrive at the following inequality, valid for $B
> 2 c_2$,
$$\pr{M_n(\eps) > B\ \eps^{(d-k)w} n } \leq \exp\left(c_1
  \eps^{-k/\alpha} - c_3\ B\ \eps^{(d-k)w} n\right).
$$

Choose $\eps = n^{-\alpha/(k+\alpha (d-k) w)}$, so that
$\eps^{-k/\alpha} = \eps^{(d-k)w} n = n^{\rho(r_0)}$.
Then, the result above transforms into 
$$\pr{M_n(\eps) > B\ n^{\rho(r_0)}} \goto 0, \quad n \goto
\infty,$$
valid for $B > (c_1/c_3) \vee 2 c_2$.

\section{Proof of Theorem \ref{th:velocity_glrt_lb}}
\label{proof:velocity_glrt_lb}
We use the notations appearing in Section
\ref{proof:velocity_glrt_ub}, except for the various constants which
are refreshed in this section.

For each $\bs \in \bS$, take $\psi_\bs:\bbR^k \goto \bbR$ infinitely
differentiable, supported in 
$[-1/2,1/2]^k$, and satisfying $\psi_\bs^{(\bt)}(0) = 0$ if $\bt \in
\bS$ and $\bt \neq \bs$; and $\psi_\bs^{(\bs)}(0) = 1$.

Let $c_1 \geq 1$ such that $c_1 \geq \|\psi_\bs^{(\bt)}\|_\infty$, for all
$\bs \in \bS$ and $\bt \in \bbN^{k}$ with $|\bt| \leq r+1$. 

Again, choose $\eps > 0$ such that $\eps^{-k/\alpha} = \eps^{(d-k)w}n$.
Define $\eps_\bs = \eps^{1-|\bs|/\alpha}$ and, with $c_2 > 1$ to be determined later, let $\eps' = (c_2 \eps)^{1/\alpha}$.
We partition $[0,1]^k$ into hypercubes of sidelength $\eps'$, 
that we denote by $I_m$, where $m = (m_1,\dots,m_k) \in ([0,1/\eps'] \cap \bbN)^k$.
For a given $m$, let 
$$R_m = I_m \times [\eps/2,\eps]^{d-k} \times \prod_{\bs \in
  \bS\setminus \{\bzero\}} [0,\eps_{\bs}]^{d-k} \subset \bbR^k \times
\bbR^{(d-k)|\bS|}.$$ 

Denote by $\cM = ([0,1/\eps'] \cap 2 \bbN)^k$ and let $M$ be the (random) subset of $m \in \cM$ such that there is $i(m) \in \{1,\dots,n\}$ with
$(X_{i(m)},Y^{\bS}_{i(m)}) \in R_m$.

\begin{lemma}
\label{lem:boundM}
There is a constant $A = A(k,d,\alpha,\beta,r_0)> 0$ such that
$$\pr{|M| > A\ n^{\rho(r_0)}} \goto 1, \quad n \goto \infty.$$ 
\end{lemma}
Lemma \ref{lem:boundM} is proved in Appendix \ref{proof:boundM}.
Theorem \ref{th:velocity_glrt_lb} now follows if we are able to interpolate the points $\{(X_{i(m)},Y^{\bS}_{i(m)}): m \in M\}$ with a function in $\cH$.\\
 
Define, for each $j=1,\dots,d-k$,
$$h_j(x) = \sum_{m \in M}\ h_{j,m}(x),$$
where 
$$h_{j,m}(x) = g_{j,m}\left(\frac{x-X_{i(m)}}{\eps'}\right),$$
with
$$g_{j,m}(x) = Y^{\bzero}_{i(m),j}\ \psi_{\bzero}(x)\ \sum_{\bs \in \bS} (\eps')^{|\bs|}\ \frac{Y^{\bs}_{i(m),j}}{Y^\bzero_{i(m),j}}\
\psi_\bs(x).$$
Note that, for $\bt \in \bbN^k$, 
$$g_{j,m}^{(\bt)}(x) = Y^{\bzero}_{i(m),j}\ \sum_{\bzero \leq {\bf t'} \leq \bt} \bino{\bt}{{\bf t'}}\  \psi_{\bzero}^{({\bf t'})}(x)\ \sum_{\bs \in \bS} (\eps')^{|\bs|}\ \frac{Y^{\bs}_{i(m),j}}{Y^\bzero_{i(m),j}}\
\psi_\bs^{(\bt - {\bf t'})}(x),$$
where if $\bs = (s_1,\dots,s_k)$ and $\bt = (t_1,\dots,t_k)$,
$$\bino{\bt}{\bs} = \bino{t_1}{s_1} \cdots \bino{t_k}{s_k}.$$
Hence, $g_{j,m}^{(\bs)}(0) = Y^{\bs}_{i(m),j}$.
This implies that $h_{j,m}^{(\bs)}(X_{i(m)}) = Y^{\bs}_{i(m),j}$, and so $h_j^{(\bs)}(X_{i(m)}) = Y^{\bs}_{i(m),j}$, since for each $j=1,\dots,d-k$, the $h_{j,m}$'s have non-overlapping supports.
Therefore, if we let $h = (h_1,\dots,h_{d-k})$, we have that $h$ interpolates each point in $\{(X_{i(m)},Y^{\bS}_{i(m)}): m \in M\}$.

Remains to prove that, for each $j=1,\dots,d-k$, $h_j \in \Hold^{k,1}(\alpha,\beta)$.  Fix $j = 1,\dots,d-k$.
Again, because for every $x \in [0,1]^k$ there is at most one $m$ such that $h_{j,m}(x) \neq 0$, we have for all $\bt \in \bbN^k$, 
$$\|h_j^{(\bt)}\|_\infty 
\leq \max_{m \in M} \|h_{j,m}^{(\bt)}\|_\infty.$$
Fix $m \in M$.
We have $\|h_{j,m}^{(\bt)}\|_\infty = \|g_{j,m}^{(\bt)}\|_\infty\ (\eps')^{-|\bt|}$,
with
$$\|g_{j,m}^{(\bt)}\|_\infty \leq Y^{\bzero}_{i(m),j}\ \sum_{\bzero \leq {\bf t'} \leq \bt} \bino{\bt}{{\bf t'}}\  \|\psi_{\bzero}^{({\bf t'})}\|_\infty\ \sum_{\bs \in \bS} (\eps')^{|\bs|}\ \frac{Y^{\bs}_{i(m),j}}{Y^\bzero_{i(m),j}}\ \|\psi_\bs^{(\bt- {\bf t'})}\|\infty.$$
Since $Y^\bzero_{i(m),j} \leq \eps$ and
$$0 \leq (\eps')^{|\bs|}\ \frac{Y^{\bs}_{i(m),j}}{Y^\bzero_{i(m),j}} \leq 2 c_2^{|\bs|/\alpha} \leq 2 c_2^{r/\alpha},$$
we have, for $|\bt| \leq r+1$, $\|g_{j,m}^{(\bt)}\|_\infty \leq c_3\ c_2^{r/\alpha}\ \eps$, with $c_3 = c_3(\alpha,k) > 0$.
In particular, $c_3$ does not depend on the choice of $c_2$.
Choose $c_2 > 1$ such that $c_3\ c_2^{r/\alpha-1} \leq \beta$.

Hence, for all $\bt$ such that $|\bt| \leq r+1$,  
$$\|h_{j,m}^{(\bt)}\|_\infty \leq c_3\ c_2^{r/\alpha}\ \eps\ (\eps')^{-|\bt|} = c_3\ c_2^{(r-|\bt|)/\alpha}\ \eps^{1-|\bt|/\alpha}.$$

This implies that, for $\eps$ small enough, $h_j$ takes values in $[0,1]$ and   
$$\|h_j^{(\bt)}\|_\infty \leq \beta,$$ 
for all $\bt$ such that $|\bt| \leq r$. 

Remains to prove that, for all $\bt$ such that $|\bt| = r$,  
$$
|h_j^{(\bt)}(x') - h_j^{(\bt)}(x)|
\leq \beta \|x' - x\|_\infty^{\alpha-r},
$$
for all $x,x' \in [0,1]^k$. 

\bitem
\item
Suppose $\|x' - x\|_\infty > \eps'$;
\beqar
|h_j^{(\bt)}(x') - h_j^{(\bt)}(x)|
& \leq & \|h_j^{(\bt)}\|_\infty\\ 
& \leq & c_3\ c_2^{r/\alpha}\ \eps (\eps')^{-r}\\ 
&\leq & c_3\ c_2^{r/\alpha}\ \eps (\eps')^{-r} (\|x - x'\|_\infty/\eps')^{\alpha - r}\\
& \leq & c_3\ c_2^{r/\alpha-1}\  \|x - x'\|_\infty^{\alpha - r}. 
\eeqar
\item
Suppose $\|x - x'\|_\infty \leq \eps'$ and let $\bt_+ = (t_1+1,\dots,t_k)$;
\beqar
|h_j^{(\bt)}(x') - h_j^{(\bt)}(x)|
& \leq & \|h_j^{(\bt_+)}\|_\infty \|x - x'\|_\infty\\
& \leq & c_3\ c_2^{r/\alpha}\ \eps (\eps')^{-(r+1)} \cdot (\eps')^{1 - (\alpha - r)}
\|x - x'\|_\infty^{\alpha - r}\\
& = & c_3\ c_2^{r/\alpha-1}\ \|x - x'\|_\infty^{\alpha - r}. 
\eeqar
\eitem
Since we chose $c_2$ such that $c_3\ c_2^{r/\alpha-1} \leq \beta$, we are finished.

\section{Proof of Theorem \ref{th:direction_glrt_ub}}
\label{proof:direction_glrt_ub}

Let $\Psi$ be the discrepancy on $\bbR^d \times \bbG(k,d)$ defined by
$$\Psi((z,H),(z_1,H_1)) = \max\{\|z-z_1\|_\infty,\ang{H}{H_1}^2\}.$$ 

Each $f \in \cF$ is identified with $(f,\vec{f})$,
viewed as a function on $[0,1]^k$ with values in $\bbR^d \times \bbG(k,d)$ -- the first and third constraints on the derivatives of $f \in \cF$ guarantee that $\vec{f}(x)$ is indeed a $k$-dimensional subspace of $\bbR^d$ for all $x \in [0,1]^k$.  
With this perspective, $\Psi$ induces a discrepancy on $\cF$. 
The proof is based on coverings of $\cF$ with respect
to that discrepancy -- still denoted by $\Psi$.

It turns out that $\Psi$ is dominated by the discrepancy $\Phi$
defined in Section \ref{proof:velocity_glrt_ub}, with $\alpha = 2$
and $r_0 = 1$.
Indeed, we have the following.

\begin{lemma}
\label{lem:phi_psi}
There is a constant $c = c(k,d,\beta)$ such that, for any $f,g \in
\cF$ and $x \in [0,1]^k$,
$$\ang{\vec{f}(x)}{\vec{g}(x)} \leq c\ \max_{s=1,\dots,k}
\|\partial_s f (x)-\partial_s g (x)\|_\infty.$$ 
\end{lemma}
To get Lemma \ref{lem:phi_psi}, we apply Lemma  \ref{lem:ang_spc} in Appendix \ref{app:ang_spc} with $u_i$ (resp. $v_i$) defined as $\partial_i f (x)$ (resp. $\partial_i g (x)$) and $c_1 = 1/\beta$, $c_2 = \beta$, $c_3 =  1/(2 \beta (d-k))$.   

Therefore, the $\eps$-covering number of $\cF$ with respect to $\Psi$
is bounded by the $\eps$-covering number of $\cF$ with respect to
$\Phi$, whose logarithm is of order $\eps^{-k/2}$ -- see Lemma \ref{lem:entropy}, where $d$ enters only in the constant.\\

Following the steps in Section
\ref{proof:velocity_glrt_ub}, we only need to find an equivalent of Lemma \ref{lem:graphS_measure}, namely compute an upper bound on the measure of the $\eps$-neighborhood of 
$$\{(f(x),\vec{f}(x)): x \in [0,1]^k\} \subset \bbR^d \times \bbG(k,d)$$ 
for the discrepancy $\Psi$, valid for all $f \in \cF$.
For $H \in \bbG(k,d)$, let
$$B(H,\eps) = \{K \in \bbG(k,d): \ang{H}{K} \leq \eps\}.$$
As in the proof of Lemma \ref{lem:graphS_measure}, we are left with computing a upper bound on $\lambda(B(H,\eps))$, which is independent of $H \in \bbG(k,d)$ since $\lambda$ is invariant under the (transitive) action of the orthogonal group.
(Remember that $\lambda$ denotes the uniform measure on $\bbG(k,d)$.)

\begin{lemma}
\label{lem:g_vol_ub}
There is a constant $c = c(k,d)$ such that, for all $\eps > 0$ and for all $H \in \bbG(k,d)$, $\lambda(B(H,\eps)) \leq c\ \eps^{(d-k)k}$.
\end{lemma}
Lemma \ref{lem:g_vol_ub} is a direct consequence of Lemma \ref{lem:g_disj} in Appendix
\ref{app:g_cover} and the fact that $\lambda\left(\bbG(k,d)\right) = 1$. 

\section{Proof of Theorem \ref{th:direction_glrt_lb}}
\label{proof:direction_glrt_lb}

We show that Theorem \ref{th:velocity_glrt_lb} implies Theorem \ref{th:direction_glrt_lb}. 
We first start by showing that functions of the form $f(x) = (x,g(x))$,  with $g \in \cH$, belong to $\cF$.

\begin{lemma}
\label{lem:HinF}
For all $g \in \cH$, the function $f(x) = (x,g(x))$
belongs to $\cF$.
\end{lemma}
Lemma \ref{lem:HinF} is proved in Appendix \ref{proof:HinF}.

Let $W$ be sampled uniformly at random in $\bbG(k,d)$.
With probability one, there is a unique set of vectors in
$\bbR^{d-k}$, $\{Y^{\bs}: |\bs| = 1\}$, such that $W =
\Span{(\bs,Y^{\bs}): |\bs| = 1}$.  
Indeed, $W$ has the same distribution as $\Span{w_1,\dots,w_k}$, where $w_1,\dots,w_k$ are i.i.d. uniformly distributed on the unit sphere of $\bbR^d$ and therefore, with probability one, $\inner{w_i}{e_i} \neq 0$ for all $i=1,\dots,k$, $\{e_i:i=1,\dots,d\}$ being the canonical basis of $\bbR^d$.
The uniqueness comes from the fact that a subspace of the form $\Span{(\bs,Y^{\bs}): |\bs| = 1}$ does not contain a vector of the form $(0,Y)$, with $Y \in \bbR^{d-k} \setminus \{0\}$, so it does not contain two distinct vectors of the form $(\bs,Y_1)$ and $(\bs,Y_2)$.  


Through the map $\kappa$ that associates $W$ to $\{Y^{\bs}: |\bs| = 1\}$, the uniform measure on $\bbG(k,d)$ induces a probability measure $\nu$ on $\bbR^{(d-k)k}$.
If we observe $\{(Z_i,W_i):i=1,\dots,n\}$, we let
$$Z_i = (X_i,Y_i) \in \bbR^k \times \bbR^{d-k} \quad {\rm and} \quad
W_i = \Span{(\bs,Y_i^{\bs}): |\bs| = 1}.$$
With here $\bS = \{\bs \in \bbN^k: |\bs| \leq 1\}$, we thus obtain $\{(X_i,Y_i^\bS):i=1,\dots,n\}$, independent and with common distribution $\Lambda_k \otimes \Lambda_{d-k} \otimes \nu \equiv \Lambda_d \otimes \nu$, where $\Lambda_\ell$ is the uniform measure on $[0,1]^\ell$.
Note that, if $g \in \cH$ interpolates $\{(X_i,Y_i^\bS):i=1,\dots,n\}$, then $f$ defined by $f(x) = (x,g(x))$ belongs to $\cF$ by Lemma \ref{lem:HinF} and interpolates $\{(Z_i,W_i):i=1,\dots,n\}$.
With $\Lambda_d \otimes \nu$ playing the role of $\mu$, the uniform measure on $[0,1]^k \times [0,1]^{d-k} \times
[-\beta,\beta]^{(d-k)(|\bS|-1)}$, the present setting parallels the situation in Section 4.
Following the arguments given there, we are only left with obtaining the equivalent of Lemma \ref{lem:boundM}.
Looking at the proof of Lemma \ref{lem:boundM} in Section \ref{proof:boundM}, all we need is a lower bound of the form
$$\Lambda_d \otimes \nu (R_\bzero) \geq c\ \mu(R_\bzero) = c\ (\eps')^k \eps^{(d-k)w} = c\ \eps^{k/2 + (d-k)(1+k/2)}.$$
(Here $\alpha = 2$ and $w = 1+ k/2$.) 
Because 
$$\Lambda_d \otimes \nu (R_\bzero) = c\ \eps^{k/2 + (d-k)} \cdot \nu\left([0,\eps^{1/2}]^{(d-k)k}\right),$$
the following lemma provides what we need.

\begin{lemma}
\label{lem:xi_lb}
There is a constant $c = c(k,d) > 0$ such that, for all $\eps > 0$ small
enough, 
$$\nu\left([0,\eps]^{(d-k)k}\right) > c\ \eps^{(d-k)k}.$$
\end{lemma}
To prove Lemma \ref{lem:xi_lb}, we first show that for some constant $c = c(k,d) > 0$, $\kappa^{-1}([0,\eps]^{(d-k)k})$ contains $B(H,c\, \eps)$, where $H = \Span{e_1,\dots,e_k}$ and $\{e_1,\dots,e_d\}$ is the canonical basis of $\bbR^d$. 
Indeed, let $c$ be the constant provided by Lemma \ref{lem:spc_ang} and consider $K \in B(H,c\, \eps)$.
Let $\kappa(K) = \{y^1,\dots,y^k\}$, so that $K = \Span{e_i + y^i:i=1,\dots,k}$.
Applying Lemma \ref{lem:spc_ang} with $u_i = e_i$ for $i=1,\dots,k$, $v = e_j + y^j$ and $u_{k+1} = (v - Pv)/\|v - Pv\|_\infty$, where $P$ denotes the orthogonal projection onto $H$, we get that $\ang{v}{H} \geq c\ \|v - Pv\|_\infty$.
With $\|v - Pv\|_\infty = \|y^j\|_\infty$ and the fact that $\ang{v}{H} \leq \ang{K}{H}$, we see that we have $\|y^j\|_\infty \leq \eps$.
This being true for all $j$, we have $\kappa(K) \in [0,\eps]^{k(d-k)}$.
We then apply the following result.

\begin{lemma}
\label{lem:g_vol_lb}
There is a constant $c = c(k,d) > 0$ such that, for all $\eps > 0$ and for all $H \in \bbG(k,d)$, $\lambda(B(H,c\, \eps)) \geq c\ \eps^{(d-k)k}$.
\end{lemma}
Lemma \ref{lem:g_vol_lb} is a direct consequence of Lemma \ref{lem:g_cover} in Appendix
\ref{app:g_cover} and again the fact that $\lambda\left(\bbG(k,d)\right) = 1$.

\begin{appendix}

\section{Proof of Lemma \ref{lem:boundM}}
\label{proof:boundM}

Lemma \ref{lem:boundM} is a conditional version of {\it The Coupon Collector's Problem} -- see e.g. \cite{Mos}.   
We nevertheless provide here an elementary proof.\\

Let $K = N_n^{(r_0)}(\cup_{m \in \cM}\ R_m)$.
We know that $K \sim \Bin(n,p)$, where $p = |\cM|\ \mu(R_\bzero)$ with $|\cM| \propto n^{\rho(r_0)}$ and
$$
\mu(R_\bzero) = (\eps')^k \cdot (\eps/2)^{d-k} \cdot \prod_{\bs \in
  \bS \setminus \{\bzero\}} (\eps_{|\bs|}/(2\beta))^{d-k} \propto
(\eps')^k \eps^{(d-k)w}.
$$
This implies that $p n = c' |\cM|$ for some $c' > 0$ not depending on $n$, by definition of $\eps$ and $\eps'$.
Let $c = c'/2$ and $c_0 \in (e^{-c},1)$, and also, to simplify notation, let $\ell = |\cM|$ and $S = |\cM| - |M|$.
Because $|\cM| \propto n^{\rho(r_0)}$, it is enough to show that $\pr{S > c_0\ \ell} \goto 0$ as $n \goto \infty$.

We have
$$\pr{S > c_0\ \ell} \leq \pr{S > c_0\ \ell|K > c\ \ell} + \pr{K > c\ \ell},$$ 
with $\pr{K > c\ \ell} \goto 0$ as $n \goto \infty$ by Lemma \ref{lem:bin_tail}, and
$$\pr{S > c_0\ \ell|K > c\ \ell} \leq \pr{S > c_0\ \ell|K = \lceil c\ \ell \rceil}.$$ 
Using Chebychev's inequality, we get
$$\pr{S > c_0\ \ell|K = \lceil c\ \ell \rceil} \leq \frac{\var{S|K=\lceil c\ \ell
    \rceil}}{(c_0\ \ell - \expect{S|K=\lceil c\ \ell \rceil})^2}.$$ 
We know that for any non-negative integer $k$,
$$\expect{S|K=k} = \ell (1 - 1/\ell)^k,$$
and
$$\var{S|K=k} = \ell ((1-1/\ell)^k -
(1-1/\ell)^{2k}) + \ell(\ell-1) ((1-2/\ell)^k - (1-1/\ell)^{2k}).$$
Therefore, when $\ell \goto \infty$,
$$\expect{S|K=\lceil c\ \ell \rceil} \backsim e^{-c} \ell,$$
and, for all $\ell$,
$$\var{S|K=\lceil c\ \ell \rceil} \leq c_1 \ell,$$
so that, when $\ell$ is large,
$$\frac{\var{S|K=\lceil c\ \ell \rceil}}{(c_0\ \ell -
  \expect{S|K=\lceil c\ \ell \rceil})^2} \leq \frac{c_2}{\ell}.$$
Since $\ell$ is an increasing function of $n$ that tends to infinity, we conclude that
$$\pr{S > c_0\ \ell} \goto 0, \quad n \goto \infty.$$ \qed

\section{Coverings of $\bbG(k,d)$}
\label{app:g_cover}

\begin{lemma}
\label{lem:g_disj}
There is a constant $c  = c(k,d) > 0$ such that, for all $\eps > 0$, there is $H_1,\dots,H_\ell \in \bbG(k,d)$ with $\ell > c\ \eps^{-(d-k)k}$ and $B(H_i,\eps) \cap B(H_j,\eps) = \emptyset$ if $i \neq j$.
\end{lemma}
{\it Proof of Lemma \ref{lem:g_disj}.}
Fix $\eps > 0$ and consider 
$$H_\bn = \Span{e_i + \eps \sum_{j=k+1}^d n_{i,j} e_j: i=1,\dots,k},$$
where $\bn = (n_{i,j}:i=1,\dots,k;j=k+1,\dots,d) \in ([0,1/\eps] \cap \bbN)^{(d-k)k}$.
As there are more than $1/2\ \eps^{-(d-k)k}$ such $H_\bn$'s,
it suffices to prove that, for some constant $c > 0$,
$$\ang{H_\bm}{H_\bn} \geq c\ \eps,$$
as soon as $\bm \neq \bn$, for that would imply that the balls $B(H_\bn,c/3\ \eps)$ are disjoint when $\bn$ runs through $([0,1/\eps] \cap \bbN)^{(d-k)k}$.

Therefore, fix $\bm \neq\bn$, both in $([0,1/\eps] \cap \bbN)^{(d-k)k}$.
For $i=1,\dots,k$, let $u_i = e_i + \eps \sum_{j=k+1}^d m_{i,j} e_j$ and $v_i = e_i + \eps \sum_{j=k+1}^d n_{i,j} e_j$, where we assume, without loss of generality, that $u_1 \neq v_1$.
Now, by definition
$$\ang{H_\bm}{H_\bn} \geq \ang{v_1}{\Span{u_1,\dots,u_k}}.$$
To proceed further, we apply Lemma \ref{lem:spc_ang} with $u_1,\dots,u_k$, $u_{k+1} = (v_1 - u_1)/\|v_1 - u_1\|_\infty$ and $v = v_1$.
It is straightforward to see that the conditions are satisfied, since in particular $v_1 = u_1 + \|v_1-u_1\|_\infty u_{k+1}$.
Hence, for a constant $c > 0$ depending only on $k,d$,
$$\ang{v_1}{\Span{u_1,\dots,u_k}} \geq c\ \|v_1 - u_1\|_\infty.$$
To conclude, note that $\|v_1 - u_1\|_\infty \geq \eps$.\qed

\begin{lemma}
\label{lem:g_cover}
There is a constant $c  = c(k,d) > 0$ such that, for all $\eps > 0$, there is $H_1,\dots,H_\ell \in \bbG(k,d)$ with $\ell < c\ \eps^{-(d-k)k}$ and $\bbG(k,d) \subset \bigcup_i B(H_i,\eps)$.
\end{lemma}
{\it Proof of Lemma \ref{lem:g_cover}.}
Let $e_1,\dots,e_d$ be the canonical basis for $\bbR^d$ and let $c_1$ be the constant given by Lemma \ref{lem:span}.
As in the proof of Lemma \ref{lem:g_disj}, define
$$H_\bn^\sigma = \Span{e_{\sigma(i)} + \eps \sum_{j=k+1}^d n_{i,j} e_{\sigma(j)}: i=1,\dots,k},$$
where $\bn = (n_{i,j}:i=1,\dots,k;j=k+1,\dots,d) \in ([0,(c_1+1)/\eps) \cap \bbN)^{(d-k)k}$ and $\sigma$ is a permutation of $\{1,\dots,d\}$.
There are no more than $c\ \eps^{-(d-k)k}$ such subsets, where 
$$c = \bino{d}{k}\ (c_1+1)^{(d-k)k}.$$
We now show that there is a constant $c = c(k,d) > 0$ such that, for all $H \subset \bbG(k,d)$, there exists such an $H_\bn^\sigma$ satisfying $\ang{H_\bn^\sigma}{H} \leq c\ \eps$.

So fix $H \subset \bbG(k,d)$.
From Lemma \ref{lem:span}, it comes that there is $\sigma$, a permutation of $\{1,\dots,d\}$, and $\xi_{i,\sigma(j)} \in [0,c_1]$ for $i=1,\dots,k$ and $j=k+1,\dots,d$, such that
$$H = \Span{e_{\sigma(i)} + \sum_{j=k+1}^d \xi_{i,\sigma(j)} e_{\sigma(j)}: i=1,\dots,k}.$$
For $i=1,\dots,k$ and $j=k+1,\dots,d$, define $n_{i,j}$ to be the entire part of $\xi_{i,\sigma(j)}/\eps$, thus obtaining $\bn = (n_{i,j}) \in ([0,(c_1+1)/\eps) \cap \bbN)^{(d-k)k}$.
Applying Lemma \ref{lem:ang_spc} with, for $i=1,\dots,k$,
$$u_i = e_{\sigma(i)} + \sum_{j=k+1}^d \xi_{i,\sigma(j)} e_{\sigma(j)},$$
and
$$v_i = e_{\sigma(i)} + \eps \sum_{j=k+1}^d n_{i,j} e_{\sigma(j)},$$
we get $\ang{H_\bn^\sigma}{H} \leq c\ \eps$ for some constant $c = c(k,d) > 0$.
\qed

\section{Proof of Lemma \ref{lem:HinF}}
\label{proof:HinF}
Only the last property defining $\cF$ is non-trivial.
Fix $f \in \cH$, and let $f(x) = (x,g(x))$.
Fix $s \in \{1,\dots,k\}$ and $x \in [0,1]^k$, and let $\bs$ be the $s^{th}$ canonical basis vector of $\bbR^k$.
We have $\partial_s f (x) = (\bs,\partial_s g (x))$. 

Let $v = v_1 + v_2$ where $v_{1} = (\bs,0)$ and $v_2 = (0,\partial_s g (x))$, and pick one vector $w \in \Span{\partial_t f (x) : t \neq s}$.  

We first show that 
$$\inner{v}{w}^2 \leq
\frac{(d-k)\beta^2}{1+(d-k)\beta^2}\ \|v\|^2 \cdot \|w\|^2.$$
Indeed, since $v_1$ is orthogonal to $w$, we have 
$\inner{v}{w} = \inner{v_2}{w}$, so that, using the
Cauchy-Schwartz inequality,
$$\inner{v}{w}^2 \leq \|v_2\|^2 \cdot \|w\|^2.$$
We have $\|v\|^2 = \|v_1\|^2 + \|v_2\|^2$.
Moreover, $\|v_1\|^2 = 1$ and $\|v_2\|^2 \leq (d-k)\beta^2$, since
$\|\partial_s f (x)\|_\infty \leq \beta$.
Conclude with
\beqar
\|v_2\|^2
& = & \left(1 - \frac{1}{\|v\|^2}\right) \|v\|^2\\
& \leq  & \left(1 - \frac{1}{1 + (d-k)\beta^2}\right) \|v\|^2.
\eeqar

Since $w$ is arbitrary in $\Span{\partial_t f (x): t \neq s}$, this shows that
$$\ang{\partial_s f (x)}{\Span{\partial_t f (x): t \neq s}} \geq \acos\left(\sqrt{\frac{(d-k)\beta^2}{1+(d-k)\beta^2}}\right).$$
Furthermore,
\beqar
\acos\left(\sqrt{\frac{(d-k)\beta^2}{1+(d-k)\beta^2}}\right) 
& = & {\rm atan}\left(\frac{1}{\beta \sqrt{d-k}}\right)\\
& \geq & \frac{1}{2 \beta \sqrt{d-k}},
\eeqar
where the last inequality comes from the fact that ${\rm atan}(y) \geq y/2$ for $y \leq \pi/2$.

\section{Auxiliary Results in Euclidean Spaces}
\label{app:ang_spc}
\begin{lemma}
\label{lem:ang_spc}
Let $c_1,c_2,c_3$ be three positive constants.
Let $u_1,\dots,u_k; v_1,\dots,v_k \in \bbR^d$ such that for all
$i=1,\dots,k$,  $c_1 \leq \|v_i\|_\infty,\|u_i\|_\infty
\leq c_2$, and, if $k \geq 2$,
$$\ang{u_i}{\Span{u_j:j\neq i}} \geq c_3;$$
$$\ang{v_i}{\Span{v_j:j\neq i}} \geq c_3.$$
Then, for a constant $c$ depending only on $k,d,c_1,c_2,c_3$,
$$\ang{\Span{u_i:i=1,\dots,k}}{\Span{v_i:i=1,\dots,k}} \leq c \max_{i=1,\dots,k} \|v_i - u_i\|_\infty.$$
\end{lemma}
{\it Proof of Lemma \ref{lem:ang_spc}.}
By multiplying the constants that appear in the Lemma by constants that depend only on $d$, we can assume that the conditions in the Lemma hold for the Euclidean norm.
Throughout, let $\eps = \max_{i=1,\dots,k} \|v_i - u_i\|_\infty$.

First assume that $u_1,\dots,u_k$ (resp. $v_1,\dots,v_k$) are
orthonormal.
Take 
$$u = \sum_i \xi_i u_i \in \Span{u_i:i=1,\dots,k},$$
of norm equal to 1.
Define
$$v = \sum_i \xi_i v_i \in \Span{v_i:i=1,\dots,k}.$$
We show that
$$\acos(|\inner{u}{v}|) = O(\eps),$$
by showing that
$$\inner{u}{v} = 1 + O(\eps^2).$$
This comes from the fact that, since $\|u\| = \|v\| = 1$,
$$\inner{u}{v} = 1 - \|u - v\|^2/2,$$
and
$$\|u - v\| \leq \sum_i |\xi_i|\ \|u_i - v_i\| \leq \sqrt{d}\ \eps.$$ 

If $u_1,\dots,u_k$ (resp. $v_1,\dots,v_k$) are not orthonormal, we
make them so.
Define $a_1' = u_1$ and $a_1 = a_1'/\|a_1'\|$, and for $i=2,\dots,k$,
define 
$$a_i' = u_i - \sum_{j=1}^{i-1} \inner{u_i}{a_{j}} a_{j},$$
and $a_i = a_i'/\|a_i'\|$. 
Similarly, define $b_1' = v_1$ and $b_1 = b_1'/\|b_1'\|$, and for
$i=2,\dots,k$, define
$$b_i' = v_i - \sum_{j=1}^{i-1} \inner{v_i}{b_{j}} b_{j},$$
and $b_i = b_i'/\|b_i'\|$.

We have, for $i=1,\dots,k$, $c_1 \sin c_3 \leq \|a_i'\| \leq c_2$.
Indeed, since $a_i'$ is the difference between $u_i$ and its orthogonal projection onto $\Span{u_1,\dots,u_{i-1}}$, it follows that
$$\|a_i'\| = \|u_i\|\ \sin \ang{u_i}{\Span{u_1,\dots,u_{i-1}}},$$
with 
$$\ang{u_i}{\Span{u_1,\dots,u_{i-1}}} \geq \ang{u_i}{\Span{u_j: j \neq i}} \geq c_3 > 0.$$
In the same way, for $i=1,\dots,k$, $c_1 \sin c_3 \leq \|b_i'\| \leq c_2$.

We also have $a_i'-b_i' = O(\eps)$.
We prove that recursively.
First, $\|a_1'-b_1'\| = \|u_1-v_1\| \leq \eps$.
Assume $a_{i-1}'-b_{i-1}' = O(\eps)$.
This implies $a_{i-1} - b_{i-1} = O(\eps)$; indeed,
\beqar
a_{i-1} - b_{i-1} 
& = & \frac{\|b_{i-1}'\| a_{i-1}' - \|a_{i-1}'\|
  b_{i-1}'}{\|a_{i-1}'\| \|b_{i-1}'\|}\\
& \leq & \frac{(\|a_{i-1}'\|+O(\eps)) a_{i-1}' - \|a_{i-1}'\|
  (a_{i-1}'+O(\eps))}{c^2}\\
& = & O(\eps).
\eeqar 
Now,
$$a_i' - b_i' = u_i - v_i - (\inner{u_i}{a_{i-1}} a_{i-1} -
\inner{v_i}{b_{i-1}} b_{i-1}),$$
with $u_i - v_i = O(\eps)$ and
$$\inner{v_i}{b_{i-1}} b_{i-1} = \inner{v_i}{a_{i-1}+O(\eps)}
(a_{i-1}+O(\eps)) = \inner{v_i}{a_{i-1}} a_{i-1} + O(\eps).$$
So that
$$\inner{u_i}{a_{i-1}} a_{i-1} -
\inner{v_i}{b_{i-1}} b_{i-1} = \inner{u_i-v_i}{a_{i-1}} a_{i-1} +
O(\eps) = O(\eps).$$
Hence, the recursion is satisfied.

We then apply the first part to $a_1,\dots,a_k$ and $b_1,\dots,b_k$.\qed

\begin{lemma}
\label{lem:spc_ang}
Fix $c_1,c_2,c_3$ three positive constants.
Let $u_1,\dots,u_{k+1} \in \bbR^d$ such that for $i=1,\dots,k+1$, $c_1 \leq \|u_i\|_\infty \leq c_2$, and
$$\ang{u_i}{\Span{u_j:j\neq i}} \geq c_3.$$
Then, there is a positive constant $c$ depending only on $k,d,c_1,c_2,c_3$ such that, for all $v = \sum_i \xi_i u_i$ with $\|v\|_\infty \leq c_2$,
$$\ang{v}{\Span{u_i:i=1,\dots,k}} \geq c\ (|\xi_{k+1}| \wedge 1).$$
\end{lemma}
{\it Proof of Lemma \ref{lem:spc_ang}.}
The proof is similar to that of Lemma \ref{lem:ang_spc} above.
Again, we may work with the Euclidean norm instead of the supnorm.

First assume that $u_1,\dots,u_{k+1}$ are orthonormal.
Take $v = \sum_i \xi_i u_i$ of norm equal to $\sqrt{\sum_i \xi_i^2} \leq c_2$.
Let $Pv = \sum_{i=1}^k \xi_i u_i$, the orthonormal projection of $v$ onto $\Span{u_1,\dots,u_k}$.
By definition,
$$\ang{v}{\Span{u_1,\dots,u_k}} = \acos\left(\frac{\inner{v}{Pv}}{\|v\| \|Pv\|}\right).$$
We then conclude with
$$\frac{\inner{v}{Pv}}{\|v\| \|Pv\|} = \sqrt{1 - \frac{\xi_{k+1}^2}{\|v\|^2}} = 1 + O(\xi_{k+1}^2).$$
 
In general, we first make $u_1,\dots,u_{k+1}$ orthonormal as we did in the proof of Lemma \ref{lem:ang_spc}, except in reverse order, meaning that $a_{k+1} = u_{k+1}/\|u_{k+1}\|$.
Since for all $v = \sum_i \xi_i u_i = \sum_i \gamma_i a_i$, $|\gamma_{k+1}| = |\xi_{k+1}| \|u_{k+1}\| \geq c_1\ |\xi_{k+1}|$, we can apply the first part to $a_1,\dots,a_{k+1}$.\qed

\begin{lemma}
\label{lem:perm_bases}
Let $e_1,\dots,e_d$ be the canonical basis of $\bbR^d$.
There is a constant $c = c(k,d) > 0$ such that, for every $u_1,\dots,u_k$, orthonormal set of vectors in $\bbR^d$, there exists a permutation $\sigma$ of $\{1,\dots,d\}$ such that
$$|\inner{u_i}{e_{\sigma(i)}}| \geq c, \quad \forall i=1,\dots,k.$$  
\end{lemma}
{\it Proof of Lemma \ref{lem:perm_bases}.}
We prove Lemma \ref{lem:perm_bases} by recursion on $k$.
For $k=1$, we may choose $c(1,d) = 1/\sqrt{d}$, since $u_1$ is of norm 1.
Suppose the result is true at $k-1$, and consider the case at $k$.
Without loss of generality, we may assume that
$$|\inner{u_i}{e_i}| \geq c(k-1,d), \quad \forall i=1,\dots,k-1.$$
We need to show that there is a constant $c_1 > 0$ and $j \in \{k,\dots,d\}$ such that
$$|\inner{u_k}{e_j}| \geq c_1.$$
As Lemma \ref{lem:perm_bases} implies Lemma \ref{lem:span}, we can use the latter at $k-1$ to get vectors $e_i + v_i$, $i=1,\dots,k-1$, with   
$v_i \in \Span{e_k,\dots,e_d},$
and $\|v_i\|_\infty \leq c_2$, such that 
$$\Span{e_1+v_1,\dots,e_{k-1}+v_{k-1}} = \Span{u_1,\dots,u_{k-1}}.$$
Let $\xi = \max_{j=k,\dots,d} |\inner{u_k}{e_j}|$.
For all $i=1,\dots,k-1$, we have
$$|\inner{u_k}{v_i}| \leq (d-k) c_2 \xi.$$
Now, since for all $i=1,\dots,k-1$, $\inner{u_k}{e_i+v_i} = 0$, we also have
$$|\inner{u_k}{e_i}| \leq (d-k) c_2 \xi.$$  
Since there is $i=1,\dots,d$ such that $|\inner{u_k}{e_i}| \geq 1/\sqrt{d}$, we must have
$$\xi \geq \frac{1}{\sqrt{d}} \vee \frac{1}{\sqrt{d} (d-k) c_2}.$$
Conclude by calling the right handside $c_3$ and letting 
$$c(k,d) = c(k-1,d) \wedge c_3.$$
\qed

\begin{lemma}
\label{lem:span}
Let $e_1,\dots,e_d$ be the canonical basis of $\bbR^d$.
There is a constant $c = c(k,d) > 0$ such that, for $u_1,\dots,u_k$ any orthonormal set of vectors in $\bbR^d$,
there exists a permutation $\sigma$ of $\{1,\dots,d\}$, such that
$$\Span{u_1,\dots,u_k} =
\Span{e_{\sigma(1)}+v_1,\dots,e_{\sigma(k)}+v_k},$$
where, for all $i=1,\dots,k$,
$$v_i \in \Span{e_\sigma(j): j=k+1,\dots,d},$$
and $\|v_i\|_\infty \leq c$.   
\end{lemma}

{\it Proof of Lemma \ref{lem:span}.}
Applying  Lemma \ref{lem:perm_bases}, there is $c_1 > 0$ and 
a permutation $\sigma$ such that
$$|\inner{u_i}{e_{\sigma(i)}}| \geq c_1, \quad \forall i=1,\dots,k.$$
Without loss of generality, suppose $\sigma = {\rm id}$.

We now triangulate the matrix with column vectors $u_1,\dots,u_k$.
In other words, we consider $\{u_i':i=1,\dots,k\}$, where $u_i'$ is the orthogonal projection of $u_i$ onto $\Span{e_i,e_{k+1},\dots,e_d}$. 
For all $i=1,\dots,k$, we have $u_i' = \xi_i e_i + w_i$, where  $|\xi_i| \geq c_1$ and $w_i \in \Span{e_{k+1},\dots,e_d}$ with $\|w_i\|_\infty \leq 1$. 
Define $v_i = w_i/\xi_i$ and conclude with the fact that 
$$\Span{u_1',\dots,u_k'} = \Span{u_1,\dots,u_k}.$$
\qed

\end{appendix}

\def\cprime{$'$} \def\lfhook#1{\setbox0=\hbox{#1}{\ooalign{\hidewidth
  \lower1.5ex\hbox{'}\hidewidth\crcr\unhbox0}}}

\end{document}